\documentclass[12pt]{amsart}
\usepackage{amsmath,amsthm,amsfonts,amscd,amssymb,eucal,latexsym,mathrsfs}
\usepackage[all]{xy}

\usepackage[english, french]{babel}

\newtheorem{theorem}{Th\'eor\`eme}

\newtheorem{corollary}[theorem]{Corollaire}
\newtheorem{lemme}[theorem]{Lemme}
\newtheorem{proposition}[theorem]{Proposition}

\theoremstyle{definition}
\newtheorem{definition}[theorem]{D\'efinition}
\newtheorem{remark}[theorem]{Remarque}

\theoremstyle{plain}

\theoremstyle{definition}

\theoremstyle{remark}

\usepackage{graphics,graphicx}

\newcommand{\supp}{{\rm supp}}

\newcommand{\IN}{{\mathbb N}}

\newcommand{\iso}{{\rm i}}

\newcommand{\SL}{{\rm SL}}

\newcommand{\G}{\Gamma}

\renewcommand{\t}{{\mathfrak{s}}}
\renewcommand{\iso}{{\mathfrak{iso}}}
\newcommand{\sat}{{\mathrm{sat}}}

\newcommand{\CI}{{\mathbb C}}\newcommand{\RI}{{\mathbb R}}
\newcommand{\ZI}{{\mathrm{\bf Z}}}\newcommand{\NI}{{\mathrm{\bf N}}}

\allowdisplaybreaks

\date{\today}
\usepackage{stmaryrd}

\title{La propri\'et\'e de d\'ecroissance rapide pour le groupe de Wise}

\author{Sylvain Barr\'e}
\author{Mika\"el Pichot}

\begin{document}

\selectlanguage{french}

\maketitle

\begin{abstract}
On montre que le groupe $G$ de pr\'esentation 
\[
\langle a,b,c,s,t\mid c=ab=ba,\, c^2=sas^{-1}=tbt^{-1}\rangle
\]
(introduit par D. Wise) a la propri\'et\'e de d\'ecroissance rapide de Haagerup--Jolissaint, et qu'il v\'erifie donc la conjecture de Baum-Connes.
\end{abstract}

Le groupe $G$ de  pr\'esentation
\[
\langle a,b,c,s,t\mid c=ab=ba,\, c^2=sas^{-1}=tbt^{-1}\rangle
\] 
a \'et\'e introduit par D. Wise  dans \cite{Wise}. Il montre que $G$ est non Hopfien, et donc non r\'esiduellement fini.  

Notons $\ell$ la longueur des mots associ\'ee \`a la pr\'esentation de $G$ donn\'ee.

On dit que $G$ a la \emph{propri\'et\'e de d\'ecroissance rapide} relativement \`a $\ell$ s'il existe un polyn\^ome  $P$ tel que, pour tout $r\in \RI_+$
et $f,g\in \CI G$  tels que $\supp(f)\subset B_r$, on a 
\[
\|f*g\|_2\leq P(r)\|f\|_2\|g\|_2
\]  
o\`u $B_r=\{x\in G,~\ell(x)\leq r\}$ est la boule de rayon $r$ dans $G$, $\supp(f)$ est l'ensemble des \'el\'ements $x\in G$ tels que $f(x)\neq 0$,  et $\CI G$ l'alg\`ebre de groupe \`a coefficients complexes de $G$. 

Nous renvoyons \`a l'article d'Alain Valette \cite{Val-bc} pour une introduction \`a la propri\'et\'e de d\'ecroissance rapide. 
Par exemple:
\begin{itemize} 
\item les groupes \`a croissance polynomiale, les groupes libres, les groupes hyperboliques, ainsi que certains r\'eseaux uniformes (essentiellement dans $\SL_3$) satisfont cette propri\'et\'e,
\item les groupes moyennables \`a croissance non polynomiale, les r\'eseaux non uniformes, tels que $\SL_3(\ZI)$, ne la satisfont pas.
\end{itemize} 

Une conjecture  de Valette affirme que la propri\'et\'e de d\'ecroissance rapide est satisfaite pour tout  r\'eseau uniforme dans un groupe de Lie semi-simple.
R\'epondant \`a une  question de Mark Sapir, nous  montrons :

\begin{theorem}\label{th1}
Le groupe $G$ a la propri\'et\'e de d\'ecroissance rapide relativement \`a $\ell$.
\end{theorem}

Le corollaire suivant est une application imm\'ediate des travaux de   Lafforgue \cite{Laf-bc}.

\begin{corollary}\label{c2}
Le groupe de Wise $G$ satisfait \`a la conjecture de Baum-Connes,   i.e.\ l'application de Baum-Connes (sans coefficient)  
\[
\mu_{r} : K_*^\mathrm{top} (G) \to K_*(C^*_{r}(G))
\]
est un isomorphisme.
\end{corollary}

La preuve du th\'eor\`eme \ref{th1} repose sur l'\'etude des propri\'et\'es ``de rang interm\'ediaire" du groupe $G$ (au sens o\`u nous l'entendons dans \cite{rd}). Il faut y ajouter un ingr\'edient suppl\'ementaire, une propri\'et\'e dite ``de prolongement analytique" (d\'ecrite dans la section  \ref{analytique} ci-dessous), qui n'est pas satisfaite en g\'en\'eral pour  les groupes \'etudi\'es dans \cite{rd}.

\begin{remark} \begin{enumerate}
\item Le pr\'esent article est une mise \`a jour d'une courte note des auteurs (non publi\'ee) portant le m\^eme titre  \cite{wise2007}, \'ecrite fin 2007. Le contenu de la nouvelle version (2012) est essentiellement identique, mais la pr\'esentation est plus d\'etaill\'ee.

\item  Le probl\`eme de d\'eterminer si $G$ poss\`ede la propri\'et\'e de d\'ecroissance rapide a \'et\'e soulev\'e par exemple dans \cite{CSC}. Comme mention\'e dans \cite{questions}, le groupe $G$ n'est pas relativement hyperbolique et il n'est pas non plus un groupe  CAT(0) cubique. Par ailleurs, n'\'etant pas r\'esiduellement fini, il n'est pas isomorphe \`a un sous groupe de type fini de  $\mathrm{SL}_3$ sur un corps (d'apr\`es un th\'eor\`eme de Malcev), et il ne satisfait donc pas aux crit\`eres connus pour la propri\'et\'e de d\'ecroissance rapide. Parfois,  le groupe $G$ \'etait propos\'e comme contre-exemple potentiel \`a cette propri\'et\'e  (voir par exemple la Question 6.6 de \cite{questions}).    
\end{enumerate}
\end{remark}

Les  sections qui suivent contiennent la preuve du th\'eor\`eme. La propri\'et\'e de d\'ecroissance rapide est obtenue aux sections \ref{S-deduire RD} et \ref{S-saturation isoc}  en appliquant les crit\`eres connus. Les quatres premi\`eres sections permettent de se ramener \`a ces crit\`eres.   Nous montrons en fait un r\'esultat plus fort, qui affirme que \emph{le groupe de Wise $G$ est \`a branchement polynomial}, au sens de \cite{rd}, d\'efinition 16 (voir le th\'eor\`eme  \ref{wise-branchpol} ci-dessous).

\section{Propri\'et\'es du link}\label{S-link}

\smallskip

La pr\'esentation  de $G$ donn\'ee ci-dessus d\'etermine un complexe $X$ \`a courbure n\'egative (au sens CAT(0))  de dimension 2 et une action libre cocompacte de $G$ avec une orbite de sommets (voir \cite{Wise}). 

On rappelle que le link en un sommet de $X$ est le graphe m\'etrique donn\'ee par la trace sur $X$ d'une sph\`ere de petit rayon centr\'ee en ce sommet, munie de la m\'etrique angulaire. 

Le link $L$ aux sommets de $X$ est repr\'esent\'e sur la figure suivante.

\begin{figure}[htbp]
\centerline{\includegraphics[width=5.5cm]{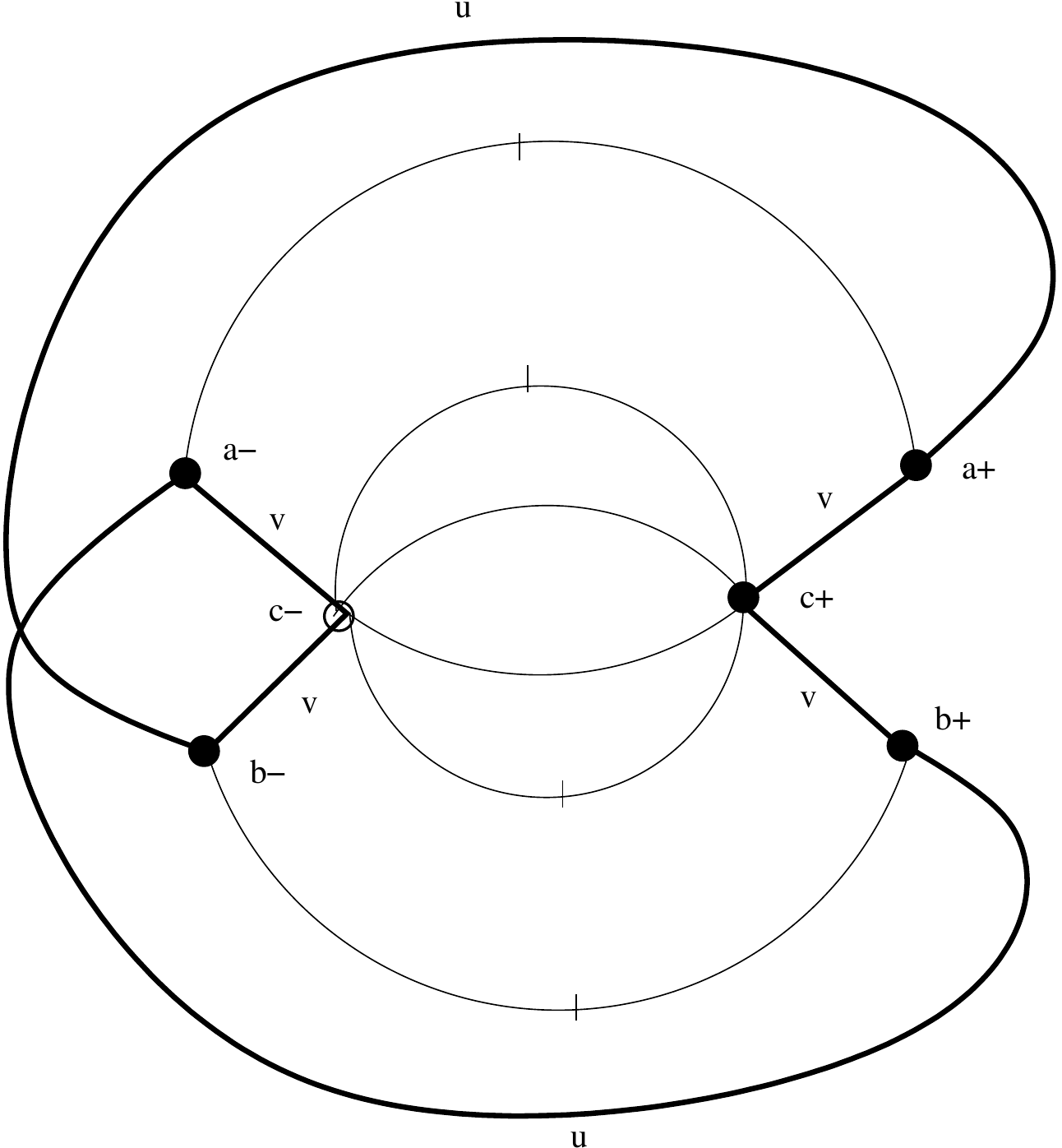}}
\caption{Le link  du groupe de Wise}\label{fig2}
\end{figure}

Le calcul du link se fait directement \`a partir de la pr\'esentation CAT(0) de $X$ telle que donn\'ee dans \cite{Wise}.

Les faces de $X$ sont compos\'ees de triangles isoc\`eles  ayant deux c\^ot\'es de longueur 1 et un de longueur $1\over 2$, et de carr\'es dont les c\^ot\'es ont longueur 1. On note $u$ et $v$ les angles aux sommets des triangles, de sorte que $u+2v=\pi$. Dans $L$ la longueur des ar\^etes en gras est indiqu\'ee sur la Figure \ref{fig2}. Les autres ar\^etes sont de longueur $\pi$. 
Les cycles (i.e.\ courbes ferm\'ees simples) de longueur $2\pi$ de $L$ sont de trois sortes :
\begin{enumerate}
\item[(1)] le cycle gras, donn\'e par le chemin $u+2v+u+2v$,
\item[(2)] les cycles mixtes, donn\'es par les chemins $u+2v+\pi$,
\item[(3)] les cycles non gras, donn\'es par les chemins centraux de la forme $\pi+\pi$.
\end{enumerate}
Le cycle gras correspond aux pavages du plan en triangles isoc\`eles dans $X$. Les ar\^etes de longueur $\pi$  viennent compl\'eter ce cycle  pour former le link tout entier. Notons que les cycles non gras  se r\'epartissent eux-m\^emes en trois sous-cat\'egories:
 le cycle de la forme $\pi+\pi$,
 les quatre cycles de la forme $(\pi/2+\pi/2)+\pi$,
et le cycle de la forme $(\pi/2+\pi/2)+(\pi/2+\pi/2)$.


Nous aurons besoin du lemme suivant.

\begin{lemme}
Soient $\alpha,\beta$ deux points de $L$ \`a distance  $> \pi$. Alors $\alpha$ et $\beta$ appartiennent \`a un  cycle de $L$ de longueur  $2\pi+2u$ ou de longueur  $2\pi+2v$. Ce cycle  est alors l'unique plus petit cycle qui  contient $\alpha$ et $\beta$.  
\end{lemme}

\begin{proof}
Il s'agit d'une v\'erification imm\'ediate.
\end{proof}

\smallskip

\section{Classification des  chromosomes}

Nous introduisons une notion de ``chromosome" pour le complexe $X$ de Wise et classifions ces sous-structures de $X$.

\begin{definition}
On appelle \emph{bande} dans le rev\^etement universel $X$ l'adh\'erence d'une composante connexe non born\'ee du compl\'ementaire du lieu singulier. 
\end{definition}

On voit que les bandes de $X$ sont de largeur 1. 

\begin{lemme}\label{chr}
Dans  $X$ \'etant donn\'ees deux bandes $B_1$ et $B_2$ il y a trois types d'intersections possibles :
\begin{enumerate}
\item[(i)] soit $B_1$ et $B_2$ ont un bord commun,
\item[(ii)] soit $B_1$ et $B_2$ ont un unique sommet  commun $x\in X$, et elles d\'efinissent dans le link en $x$ un cycle de longueur $2\pi+2u$, 
\item[(iii)] soit $B_1$ et $B_2$ ont un unique sommet commun $x\in X$, et elles d\'efinissent dans le link en $x$ un cycle de longueur $2\pi+2v$, 
\end{enumerate}  
\end{lemme}

\begin{proof}
Soient deux bandes $B_1$ et $B_2$ et $x$ un point commun \`a $B_1$ et  $B_2$. Deux cas se pr\'esentent. 

Supposons dans un premier temps que $x$ n'est pas l'unique point commun. Notons $I$ le plus grand segment (\'eventuellement infini) contenant $x$ et inclus dans $B_1\cap B_2$. Par hypoth\`ese, ce segment n'est pas r\'eduit \`a $x$.  Il s'agit de montrer que $I$ est bi-infini. Si ce n'est pas le cas, $I$ a une  extremit\'e  $y$ dans $X$. Il est facile de voir que cette extr\'emit\'e est un sommet de $X$. De plus $y$ est soit un sommet de carr\'e, soit un milieu de carr\'e dans $B_1$ et $B_2$. Dans les deux cas, la trace de $B_1\cup B_2$ sur le link $L_y$ de $y$ est un chemin non gras connexe de longueur $2\pi$. Or, un tel  chemin est n\'ecessairement un \emph{cycle} non gras de longueur $2\pi$ dans $L_y$, ce qui contredit la maximalit\'e de $I$. Donc $I=B_1\cap B_2$ est bi-infini et nous obtenons l'assertion (i).

Supposons maintenant que $B_1\cap B_2=\{x\}$.  Notons $I_1$ et $I_2$ les deux segments de longueur $\pi$ dans $L_x$ (le link de $x$) correspondant \`a $B_1$ et $B_2$.   Alors: soit $I_1$ et $I_2$ sont dans un cycle mixte de longueur $2\pi+2u$, soit ils sont dans un cycle mixte de longueur $2\pi+2v$. Ceci correspond aux deux cas $(ii)$ et $(iii)$ ci-dessus.
\end{proof}

\begin{definition}
On appelle \emph{chromosome} de $X$ la r\'eunion de deux bandes d'intersection non vide. Un chromosome est dit respectivement \emph{coll\'e}, \emph{de type $u$}, ou \emph{de type $v$},  suivant les cas  (i), (ii) et (iii) du Lemme \ref{chr} ci-dessus.  Dans les cas (ii) et (iii)  le point $x$ est appel\'e \emph{centrom\`ere}.
 \end{definition}
 
 Les trois types de chromosomes de $X$ sont repr\'esent\'es sur la figure suivante (le troisi\`eme est de type $u$, les deux derniers de type $v$).

\begin{figure}[htbp]
\centerline{\includegraphics[width=15cm]{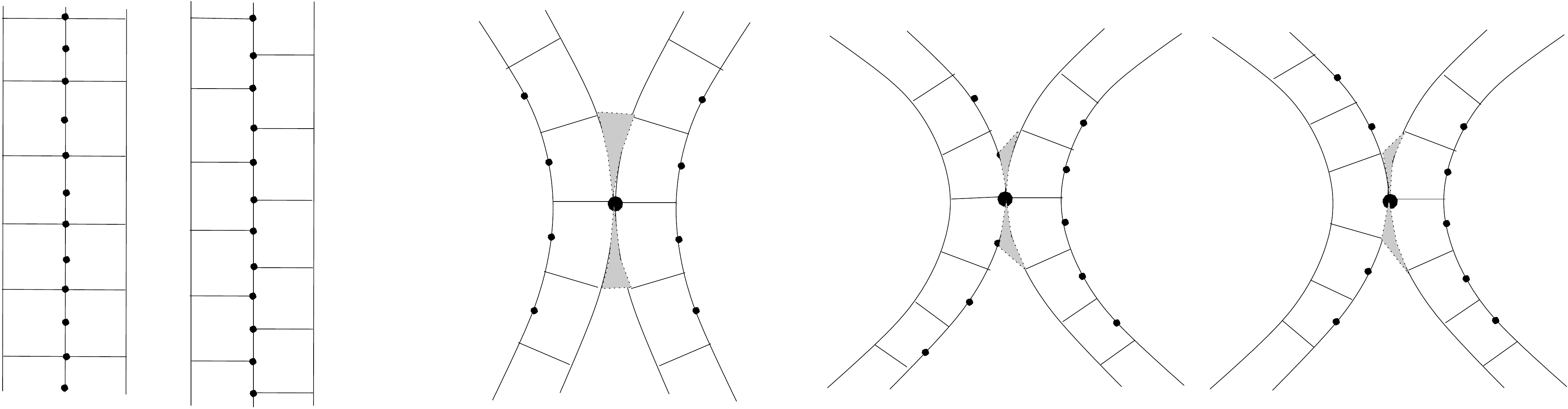}}
\caption{Chromosomes coll\'es, de type $u$,  de type $v$}\label{fig3}
\end{figure}

Soit $\gamma$ un segment g\'eod\'esique de $X$ et soit $B$ une bande. On dit que $\gamma$ \emph{rencontre} $B$ s'il contient au moins un point int\'erieur \`a $B$. On dit que $\gamma$ \emph{rencontre successivement} deux bandes $B_1$ et $B_2$  de $X$  si $\gamma$ rencontre $B_1$ et $B_2$  et si l'int\'erieur de toute autre bande $B$ de $X$,  $B\neq B_1$, $B\neq B_2$, est disjoint du sous-segment g\'eod\'esique $\gamma_0$ de $\gamma$  compris entre $B_1$ et $B_2$ (i.e.\ $\gamma_0$ ne rencontre pas $B$). 

\begin{lemme}\label{parallele}
\begin{enumerate}
\item[(i)] Soit $\Pi$ un plat de triangles isoc\`eles de $X$. Une bande dans $X$ dont l'intersection  avec $\Pi$ est non vide, et non r\'eduite \`a un point,  contient une g\'eod\'esique singuli\`ere de $\Pi$. 
\item[(ii)] Soit $\gamma$ un segment g\'eod\'esique de $X$. Si $\gamma$ ne rencontre aucune bande de $X$, alors $\gamma$ est contenu dans un  unique plat en triangles isoc\`eles  de $X$.
\item[(iii)] Soit $\gamma$ un segment g\'eod\'esique et  soient $B_1$, $B_2$ deux bandes de $X$ d'intersection  vide. Si $\gamma$  rencontre successivement $B_1$ et $B_2$,   alors il existe un plat de triangles isoc\`eles $\Pi$ dans $X$ tel que les intersections $\Pi\cap B_1$ et $\Pi\cap B_2$ sont deux droites singuli\`eres parall\`eles de $\Pi$.
\end{enumerate}
\end{lemme}

\begin{proof}
(i) Soit $B$ une bande et $\Pi$ un plat en triangles isoc\`eles tel que $I=B\cap \Pi$ contient au moins deux points. Par convexit\'e, $I$ est un segment du bord de $B$. Il suffit de montrer que $I$ est bi-infini. Sinon, notons $x$ une extr\'emit\'e de $I$ dans $X$. Alors $x$ est un sommet de $X$. De plus, nous avons un cycle de longueur $2\pi=2(u+2v)$ correspondant aux triangles isoc\`eles et un chemin de longueur $\pi$ d'origine sur ce cycle et correspondant \`a $B$. La section \ref{S-link} montre que la singularit\'e obtenue est de type livre, ce qui contredit le fait que $x$ est extr\'emal dans $B\cap \Pi$.

(ii) Si $\gamma$ est un segment g\'eod\'esique qui ne rencontre aucune bande, alors il est enti\`erement contenu dans une r\'eunion de triangles isoc\`eles.  Il suffit donc de montrer que tout point d'un triangle isoc\`ele est contenu dans un unique plat en triangle isoc\`ele. Soit $t$ un triangle isoc\`ele contenant un tel point, et soit $x$ un sommet de $t$. Alors $t$ d\'etermine une ar\^ete grasse du link $L_x$, qui se compl\`ete d'une unique fa\c con en un cycle d'ar\^etes grasses. Ce cycle gras d\'etermine un disque simplicial plat $D_1$ dans $X$ de centre $x$. On construit ensuite par r\'ecurrence une suite emboit\'ee $D_1\subset D_2\subset D_3\subset \cdots$ de disques simpliciaux plats en triangles isoc\`eles de rayons (simpliciaux) respectifs $1, 2, 3, \ldots$, de la fa\c con suivante. Supposons $D_i$ construit pour $i\geq 1$. Notons $x_i^1,\ldots x_i^{k_i}$ les sommets du bord de $D_i$ organis\'es par ordre cyclique. Pour chaque entier $j\in [1,k_i]$, la trace de $D_i$ dans $L_{x_i^j}$ est un chemin gras de longueur au plus $\pi$. Celui se compl\`ete d'exactement une fa\c con en un cycle gras de longueur $2\pi$. En outre, il est facile de voir que ces compl\'etions sont deux-\`a-deux compatibles.  Ceci montre qu'il existe exactement une extension de $D_i$ en un disque plat $D_{i+1}$ en triangles isoc\`eles de rayon simplicial $i+1$.  Par suite, $\Pi=\bigcup_{i=1}^\infty D_i$ est l'unique plat en triangles isoc\`ele contenant le point donn\'e.

(iii) Notons $C$ l'enveloppe convexe de $B_1$ et $B_2$ priv\'ee de $B_1$ et $B_2$. Comme $B_1\cap B_2=\emptyset$, $C\neq \emptyset$. D'apr\`es (ii), $C$ est inclus dans un (unique) plat en triangle isoc\`eles. Il est facile d'en d\'eduire que $B_1\cap \bar C$ et $B_2\cap \bar C$ sont deux droites parrall\`eles de ce plat.
\end{proof}

\section{Enveloppe analytique et enveloppe analytique r\'eduite d'un segment}\label{analytique}

Soit $\gamma$ un segment g\'eod\'esique de $X$. On appelle \emph{enveloppe analytique} de $\gamma$ la surface $H(\gamma)$ de $X$ (\'eventuellement \`a bords, branch\'ee, et non n\'ec\'essairement totalement g\'eod\'esique)  construite de la fa\c con suivante. 

\textbullet\ Si 
 $\gamma$ ne rencontre aucune bande de $X$, alors  $H(\gamma)$ est soit : 
 \begin{enumerate}
 \item[(1)] l'enveloppe convexe simpliciale de $\gamma$ dans le plat $\Pi$ en triangles isoc\`eles donn\'e par le lemme \ref{parallele} (ii) si $\gamma$ est non singulier, 
 \item[(2)] la r\'eunion des chromosomes coll\'es de $X$ qui intersectent $\gamma$ en au moins deux points si $\gamma$ est singulier. 
 \end{enumerate}
 
\textbullet\ Si $\gamma$ rencontre au moins une bande de $X$, notons $B_1,\ldots B_n$ les bandes de $X$ que $\gamma$ rencontrent successivement en partant de son origine  (pour une orientation fix\'ee de $\gamma$)  et proc\'edons de la fa\c con suivante. Soit $H_0(\gamma)$ la r\'eunion  des bandes $B_i$, $i=1\ldots n$. Construisons par r\'ecurrence une suite emboit\'ee de surfaces \`a bords 
\[
H_0(\gamma)\subset H_1(\gamma)\subset \ldots \subset H_{n-1}(\gamma).
\]
Supposons $H_{i-1}(\gamma)$ construite pour un entier $i\geqslant 1$ et construisons $H_i(\gamma)$ pour  $i<n$. On distingue les cas suivant.
 
\begin{enumerate}
\item[(1)] $B_i$ et $B_{i+1}$ ne s'intersectent pas : on applique la partie (iii) du Lemme \ref{parallele} en posant $H_i(\gamma)=H_{i-1}(\gamma)\cup P$, o\`u $P$ est la partie du plat $\Pi$ entre les deux droites singuli\`eres d\'efinies  dans (iii).
\item[(2)] $B_i$ et $B_{i+1}$  ont un bord commun : on pose $H_i(\gamma)=H_{i-1}(\gamma)$.
\item[(3)] $B_i$ et $B_{i+1}$ ont exactement un point commun $x\in X$ : dans ce cas $B_i\cup B_{i+1}$ est un chromosome $C$ non coll\'e  de $X$ de centrom\`ere $x$. Par construction $\gamma$ et $C$ sont transverses, i.e.\  $\gamma$ a une intersection non vide avec l'int\'erieur des deux branches de $C$, et  on compl\`ete en posant $H_i(\gamma)=H_{i-1}(\gamma)\cup P$, o\`u $P$ est la r\'eunion des deux secteurs de triangles  d'angles   $u$ (ou $v$) en $x$, d\'etermin\'es par le type $u$ (ou $v$) du chromosome $C$.
\end{enumerate}

D'apr\`es le lemme \ref{chr} il n'y a pas d'autre possibilit\'e d'intersection.

La construction ci-dessus d\'etermine  une surface $H_{n-1}(\gamma)$ dans $X$. On pose alors $H(\gamma)=H_{n-1}(\gamma)\cup P_1\cup P_2\cup Q$ o\`u $Q$ est la r\'eunion des chromosomes coll\'es de $X$ qui intersectent $\gamma$ en au moins deux points, et $P_1$, $P_2$ sont d\'efinis comme suit : si l'origine de $\gamma$ (resp. l'extr\'emit\'e de $\gamma$) n'est pas incluse dans $H_{n-1}$, alors $P_1$ (resp. $P_2$) est le demi-plan en triangles isoc\`eles bordant $B_1$ et contenant l'origine (resp. l'extr\'emit\'e) de $\gamma$ ; sinon on pose $P_1=\emptyset$ (resp. $P_2=\emptyset$).

Ceci termine la construction de $H(\gamma)$. On a choisi le terme `analytique' en r\'ef\'erence au prolongement analytique usuel: nous voyons l'enveloppe $H(\gamma)$ comme le  ``prolongement analytique" du segment $\gamma$ dans $X$. 

Enfin nous tronquons cette enveloppe :

\begin{definition}
Soit $\gamma$ un segment g\'eod\'esique de $X$. On appelle \emph{enveloppe analytique r\'eduite} de $\gamma$, et on note $H'(\gamma)$, l'ensemble des points de l'enveloppe analytique $H(\gamma)$ de $\gamma$ qui sont \`a distance au plus $2|\gamma|$ du milieu de $\gamma$, o\`u $|\gamma|$ est la longueur CAT(0) de $\gamma$ dans $X$. 
\end{definition}

\section{R\'eduction des triangles g\'eod\'esiques}

\begin{lemme}\label{reduc}
Soient $A,B,C$ trois sommets de $X$ et soit $D$ le disque totalement g\'eod\'esique de $X$ de bord le triangle g\'eod\'esique $(ABC)$. Alors l'une des deux possibilit\'es suivantes est satisfaite:
\begin{enumerate}
\item L'intersection 
\[
H'([AB])\cap  H'([BC])\cap  H'([AC])
\]
est non vide.
\item L'ensemble 
\[
D_0=D\backslash \left( H'([AB])\cup  H'([BC])\cup  H'([AC]) \right )
\]
est un disque non vide inclus dans un plat $\Pi$ en triangles isoc\`eles, dont le bord $T=\partial D_0$ est un triangle isoc\`ele g\'eod\'esique simplicial de $\Pi$. 
\end{enumerate}
\end{lemme}

Dans la plupart des cas $D_0$ est  vide, et c'est la situation (1) qui pr\'evaut.  Le lemme \ref{reduc} se d\'emontre en \'etudiant les positions possibles des chromosomes \`a l'int\'erieur de $D$. On a repr\'esent\'e un cas significatif sur la figure \ref{fig3} (o\`u les centrom\`eres de $D$ se trouvent sur les bords).

\begin{figure}[htbp]
\centerline{\includegraphics[width=12cm]{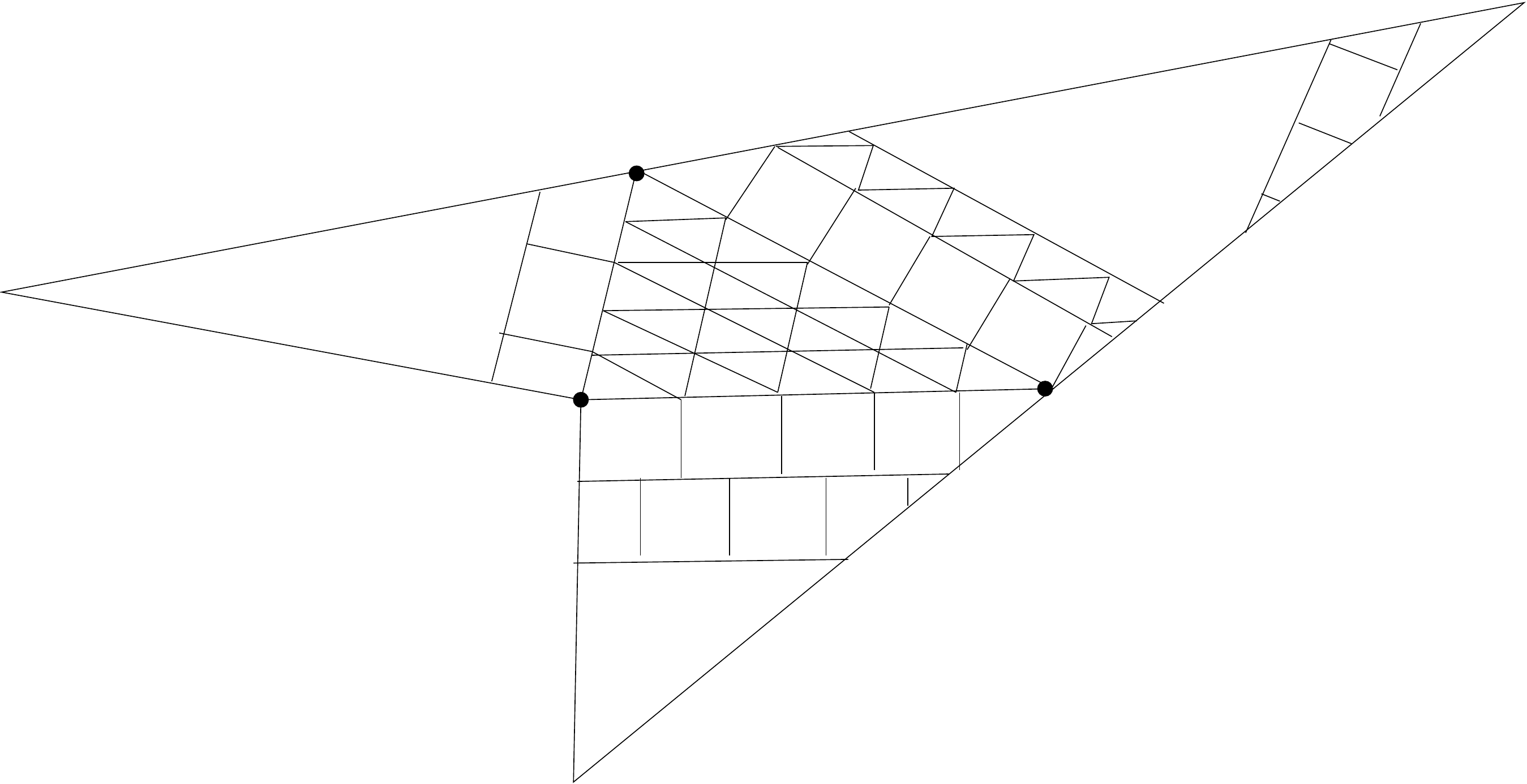}}
\caption{R\'eduction le long des enveloppes analytiques r\'eduites}\label{fig3}
\end{figure}

Avant d'aborder la d\'emonstration du lemme \ref{reduc}, nous montrons deux r\'esultats pr\'eliminaires. 
Le premier est \'el\'ementaire et g\'en\'eral:

\begin{lemme}\label{L - CAT(0) triangles }
Soient $A,B,C$ trois sommets de $X$ et soit $(ABC)$ le  triangle g\'eod\'esique correspondant. Alors les boules  centr\'ees aux milieux des cot\'es de $(ABC)$ et de rayons respectifs le double de la longueur du cot\'e correspondant \`a leur centre, ont une intersection commune contenant le petit cot\'e de $(ABC)$.
\end{lemme}

Le r\'esultat suivant est propre au complexe $X$:

\begin{lemme}\label{frizes}
Soient $A,B,C$ trois sommets de $X$. Si $[A,B]$ rencontre une bande de $X$, alors soit $[B,C]$ soit $[A,C]$ rencontre cette m\^eme bande de $X$. 
\end{lemme}

\begin{proof}
Supposons qu'il existe une bande que rencontre $[A,B]$ mais que ne rencontrent ni $[B,C]$, ni $[A,C]$. Alors les  deux points $x$ et $y$ de son bord les plus \'eloign\'es de $[A,B]$ qui appartienent au disque de bord le triangle $(ABC)$ sont int\'erieurs \`a ce disque $D$. Il en r\'esulte que le segment $[x,y]$ est singulier.  Mais cela contredit la d\'efinition d'une bande comme compl\'ementaire du lieux singulier. 
\end{proof}

\begin{proof}[D\'emonstration du lemme \ref{reduc}]
Nous distinguons plusieurs cas et sous-cas. 

(1) Supposons d'abord que l'un des c\^ot\'es, disons $[A,B]$ de $ABC$ rencontre une bande de $X$. D'apr\`es le lemme \ref{frizes},  nous pouvons supposer, quitte \`a permuter les lettres, que $[A,C]$ rencontre cette bande \'egalement. Notons $B_1$ la bande de $X$ la plus \'eloign\'ee de $A$ et qui rencontre \`a la fois $[A,B]$ et $[A,C]$, et $x$ (resp.\ $y$) le sommet de $[A,B]\cap B_1$ (resp.\ $[A,C]\cap B_1$) le plus \'eloign\'e de $A$. 

Deux cas se pr\'esentent: 
\begin{enumerate}
\item[(1a)] $[A,B]$ rencontre successivement  $B_1$ et une autre bande $B_2$ plus proche de $B$ que $B_1$. 
\item[(1b)] $[x,B]$ ne rencontre aucune bande de $X$.
\end{enumerate}  
Traitons d'abord le cas (1a). Par maximalit\'e de $B_1$ et d'apr\`es le lemme \ref{frizes}, la bande $B_2$ rencontre n\'ecessairement $[B,C]$. 

Alors:

(1a$\alpha$) Si $B_1$ et $B_2$ sont disjointes, alors elles intersectent un plat $\Pi$ en triangles isoc\`eles suivant deux droites parall\`eles (cf.\ lemme \ref{parallele} iii). Par construction, l'enveloppe convexe $\Pi'$ de ces deux droites dans $\Pi$ est incluse dans $H([A,B])$. De plus, il est facile de voir que $\Pi'$ contient aussi $[y,C]$, et qu'elle est donc incluse dans $H([A,C])$. Par suite $H([A,B])\cap H([A,C])$ contient $\Pi'$. Un raisonnement analogue montre que $H([B,A])\cap H([B,C])$ contient \'egalement $\Pi'$. Nous en d\'eduisons que $\Pi'\subset H([A,B])\cap H([A,C])\cap H([B,C])$.
Il r\'esulte alors du lemme \ref{L - CAT(0) triangles } que $H'([A,B])\cap H'([A,C])\cap H'([B,C])$ est non vide.

(1a$\beta$) Si $B_1\cap B_2\neq \emptyset$, nous appliquons la classification des chromosomes (lemme \ref{chr}) et distinguons deux sous-cas:

(1a$\beta$i) Si $B_1$ et $B_2$ forment un chromosome coll\'e,  alors $B$ appartient \`a $d=B_1\cap B_2$. Un raisonnement analogue au cas (1a$\alpha$) permet de d\'eduire que $d\subset H([A,B])\cap H([A,C])\cap H([B,C])$ et il r\'esulte du lemme \ref{L - CAT(0) triangles } que $H'([A,B])\cap H'([A,C])\cap H'([B,C])$ est non vide.

(1a$\beta$ii) Sinon, $B_1$ et $B_2$ forment un chromosome $\chi$ non coll\'e. Notons $z$ le point le plus \'eloign\'e de $B$ qui appartient \`a $[B,C]\cap B_2$. Alors $[y,z]$ est dans l'enveloppe convexe de $\chi$ et le triangle $(xyz)$ est inclus dans l'un des deux c\^ones en triangles isoc\`eles de $\chi$ issus de son centrom\`ere. Deux nouveaux cas se pr\'esentent.

(1a$\beta$ii') $[A,C]$ recontre successivement $B_1$ et une autre bande de $X$, disons $B_3$. Par maximalit\'e de $B_1$, le segment $[B,C]$ doit rencontrer successivement $B_2$ et $B_3$. De plus, les bandes $B_1$ et $B_3$ (resp.\  les  bandes $B_2$ et $B_3$) d\'eterminent un chromosome non coll\'e $\chi'$ (resp.\ $\chi''$) de $X$.  On en d\'eduit alors que le triangle $(xyz)$ est inclus dans deux c\^ones en triangles isoc\`eles de $\chi'$ et de $\chi''$ issus de leurs centrom\`eres. Donc $(xyz)$ est inclus dans $H([A,B])\cap H([A,C])\cap H([B,C])$ et il r\'esulte du lemme \ref{L - CAT(0) triangles } que $H'([A,B])\cap H'([A,C])\cap H'([B,C])$ est non vide.

(1a$\beta$ii'') Si (1a$\beta$ii') n'est pas satisfait, alors les segments $[y,C]$ et $[z,C]$ ne rencontrent pas de bande de $X$. Par construction de nos enveloppes convexes, $H([A,B])\cap H([A,C])\cap H([B,C])$ contient $(xyz)$  et il r\'esulte encore du lemme \ref{L - CAT(0) triangles } que $H'([A,B])\cap H'([A,C])\cap H'([B,C])$ est non vide.

Ceci termine la preuve dans le cas (1a). 

Passons au cas (1b). Nous avons encore deux options:

(1b$\alpha$) Supposons que $[A,C]$ rencontre sucessivement $B_1$ et une autre bande $B_2$ de $X$. D'apr\`es le lemme \ref{frizes}, $[B,C]$ rencontre $B_2$ \'egalement. Ce cas est donc un sym\'etrique du cas (1a$\beta$ii'') et se traite de la m\^eme fa\c con.

(1b$\beta$) Sinon, $[y,C]$ ne recontre aucune bande de $X$. Il en r\'esulte que $[B,C]$ ne rencontre aucune bande de $X$. Dans ce cas, le segment $[B,C]$ est inclus dans $H([A,B])\cap H([A,C])$ et le lemme \ref{L - CAT(0) triangles } permet encore de conclure que $H'([A,B])\cap H'([A,C])\cap H'([B,C])$ est non vide.

Ceci termine la preuve dans le cas (1). 

Dans le cas (2), aucun des segments $[A,B]$, $[B,C]$ et $[A,C]$ ne rencontre de bande de $X$. Il en r\'esulte que le triangle $(ABC)$ est inclus dans un plat en triangles isoc\`eles $\Pi$. Dans ce cas, on montre directement que 
\[
D_0=D\backslash \left( H'([AB])\cup  H'([BC])\cup  H'([AC]) \right )
\] 
est un triangle isoc\`ele (\'eventuellement vide) g\'eod\'esique simplicial de $\Pi$. (Une autre solution, sensiblement moins \'economique, est de d\'eformer les triangles isoc\`eles de $\Pi$  en des triangles \'equilat\'eraux et d'appliquer le lemme 29 de \cite{rd}.) 
\end{proof}

\section{En d\'eduire la propri\'et\'e de d\'ecroissance rapide}\label{S-deduire RD}

Ce qui pr\'ec\`ede permet de montrer que  $G$ v\'erifie la propri\'et\'e dite de ``branchement  polynomial" relativement \`a la m\'etrique $\ell$, au sens de la d\'efinition 16  de \cite{rd} (rappel\'ee ci-dessous). La propri\'et\'e de d\'ecroissance rapide (th\'eor\`eme \ref{th1}) r\'esultera alors directement de la proposition 17 de \cite{rd}.

Le r\'esultat principal de cet article est le suivant:

\begin{theorem}\label{wise-branchpol}
Le groupe de Wise $G$ est \`a branchement polynomial relativement \`a la m\'etrique des mots. 
\end{theorem}

Commen\c cons par des rappels terminologiques. 
Soit $\G$ un groupe d\'enombrable. Un \emph{3-chemin} de l'identit\'e $e$ \`a $z\in \G$ est un triplet  $\gamma=(a_3,a_2,a_1)$ dans $\G^3$ tel que $z=a_3a_2a_1$.  Un  triplet $(x,y,z)\in \G^3$ tel que $xy=z$ est dit \emph{triangle} de $\G$. Notons $|\cdot|$ la longueur sur $\G$.

\begin{definition}\label{HK}   
Une famille  $\G$-index\'ee de 3-chemins dans $\G$, disons $C=(C_z^r)_{z\in \G,~~r\in \NI^*}$  o\`u  $C_z^r$ un ensemble de 3-chemins de $e$ \`a $z$ dans $\G$ pour tout $z\in \G$ et $r\in \RI^*$, est dite \`a croissance polynomiale s'il existe un polyn\^ome  $p_1$ tel que pour tout $r\in \RI_+$ et tout $z\in \G$ on a 
$\#C_z^r\leq p_1(r)$.
\end{definition}

Soient $\t$ et $\t^-$ deux ensembles de triangles de $\G$ et $C=(C_z^r)_{z\in \G,~~r\in \RI_+}$ une famille $\G$-index\'ee de 3-chemins. Pour $(u,v,w)\in \t^-$ et $r\in \RI_+$ d\'efinissons $D^r_{(u,v,w)}$ comme l'ensemble des triplets $(a,b,c)$ de $\G^3$ tels que $(b^{-1},u,a)\in C^r_{b^{-1}ua}$, $(c^{-1},v,b)\in C^r_{c^{-1}vb}$ et $(c^{-1},w,a)\in C^r_{c^{-1}wa}$.

\begin{definition} \label{retr-tr}
On dit que   \emph{$\t^-$  est un retract de $\t$ le long de $C$} s'il existe un polyn\^ome $p_2$ tel que pour tout   $(x,y,z)\in \t$ il exists $(u,v,w)\in \t^-$ avec $|u|\leq p_2(|x|)$ et   $(a,b,c)\in D^{|x|}_{(u,v,w)}$ tel que $b^{-1}ua=x$ and $c^{-1}wa=z$.
\end{definition}

Rappelons la d\'efinition 16 de \cite{rd}. Un 3-chemin $(a_3,a_2,a_1)$ dans $\G$ est dit $(\kappa,\delta)$-g\'eod\'esique si $|a_1|+|a_2|+|a_3|\leq \kappa|a_3a_2a_1|+\delta$, o\`u $\kappa \geq 1$, $\delta\geq 0$ sont donn\'es.

\begin{definition}\label{polrk} On dit que $\G$ est \`a \emph{branchement polynomial} relativement \`a $|\cdot|$ s'il existe $\kappa\geq 1$, $\delta\geq 0$, une famille  $C=(C_z^r)_{z\in \G,~~r\in  \IN^*}$ d'ensembles $C_z^r$ de 3-chemins $(\kappa,\delta)$-g\'eod\'esiques de $e$ \`a $z$ \`a croissance polynomiale,  un ensemble $\t$ de triangles dans $\G$ qui sont des r\'etracts le long de $C$ de l'ensemble de tous les triangles de $\G$, et un  polyn\^ome $p_3$ tel que pour tout $z$ dans $\G$ et tout $r\in \RI_+$, le nombre de triangles dans $\t$ de la forme $(x,y,z)$ avec $|x|\leq r$ est au plus $p_3(r)$ (on dit alors que $\t$ est \`a croissance polynomiale). 
\end{definition}

On a alors: 

\begin{proposition}[Voir la prop. 17 de \cite{rd}]
Soit $\G$ un groupe d\'enombrable et $|\cdot|$ une longueur sur $\G$. Si $\G$ est \`a branchement polynomial relativement \`a $|\cdot|$, alors $\G$ satisfait \`a la propri\'et\'e de d\'ecroissace rapide relativement \`a $|\cdot|$.
\end{proposition}

Nous montrons maintenant que le groupe de Wise $G$ est \`a branchement polynomial.

Soit $z\in G$ et $r\geq 0$. Voyons $G$ comme l'ensemble des sommets de $X$ et notons $H'[z]$ l'enveloppe convexe r\'eduite du segment g\'eod\'esique de $e$ \`a $z$ dans $X$. Nous consid\'erons la famille $C_z^r$ des 3-chemins (8,0)-g\'eod\'esiques $(a_3,a_2,a_1)$  de $e$ \`a $z$ tels que $a_1,a_2a_1\in H'[z]$ et, si $|z|\geq r$, tels que $|a_1|\leq 3r$ et $|a_2|\leq r$. 

Nous avons alors:

\begin{lemme}\label{l-C polyn}
La famille $(C_z^r)_{z,r}$ est \`a croissance polynomiale. 
\end{lemme}
  
\begin{proof}  
Nous montrons qu'il existe une constante $K>0$ tel que, pour tout segment g\'eod\'esique  $\gamma$ de $X$ et tout $r>0$, le nombre de sommets de  $H'(\gamma)$ \`a distance au plus $r$ est de l'origine de $\gamma$ est major\'e par $Kr^3$. Ceci entra\^\i ne que le cardinal de $C_z^r$ est au major\'e par un polyn\^ome (de degr\'e 6), ce qu'il faut d\'emontrer.

Si $\gamma$ ne rencontre pas de bande de $X$, alors il est inclus dans un plat en triangles isoc\`eles de $X$ et le r\'esultat est clair (le polyn\^ome cherch\'e peut \^etre choisi de degr\'e 2). 

Sinon, soit $B_1, \ldots, B_n$ la suite des bandes que rencontre (successivement) $\gamma$  et qui sont \`a distance au plus $r$ de l'origine de $\gamma$. Comme les bandes sont de largeur 1, $n\leqslant r$. De plus il existe une constante $K'$ de sorte que l'ensemble $B_i^r$ des sommets de $B_i$ \`a distance au plus $r$ de $\gamma$ soit born\'e par $K'r$.   Notons $C_i$ l'ensemble des sommets de de $H(\gamma)$ situ\'es entre les bandes $B_i$ et $B_{i+1}$. D\'efinissons \'egalement $C_0$ et $C_n$ les plats en triangles contenant l'origine de $\gamma$ et, respectivement : son extr\'emit\'e si $|\gamma|\leq r$, ou bien le point de $\gamma$ \`a distance $r$ de son origine si $|\gamma|\geq r$. Comme $C_i$ est inclus dans un  plat en triangles pour tout $i=0\ldots n$, il est facile de voir qu'il existe une constante $K''$ tel que l'ensemble $C_i^r$ des sommets de $C_i$ \`a distance au plus $r$ de l'extr\'emit\'e de $\gamma$ soit de cardinal au plus $K''r^2$.   Par suite, l'ensemble des points de $H'(\gamma)$ \`a distance au plus $r$ de l'origine est born\'e par 
\[
|C_0|+ \sum_{i=1}^n |B_i^r|+|C_i^r| \leq K''r^2+(K'r+K''r^2){n(n+1)\over 2}\leq K r^3  
\]
o\`u $K$ ne d\'epend que de $K'$ et $K''$. Ceci prouve le lemme \ref{l-C polyn}.
  \end{proof}

\begin{lemme}\label{l-isopolyn}
La famille $\iso$ des triangles simpliciaux g\'eod\'esiques inclus dans un plat en triangles isoc\`eles de $X$ est \`a croissance polynomiale. 
\end{lemme}

\begin{proof}
Le lemme \ref{parallele} (ii) montre que nous pouvons en fait choisir $p_3$ constant.
\end{proof}

Le th\'eor\`eme \ref{wise-branchpol} r\'esulte alors imm\'ediatement des lemmes \ref{l-C polyn}, \ref{l-isopolyn} et du lemme  suivant:

\begin{lemme}
La famille $\iso$ des triangles simpliciaux g\'eod\'esiques inclus dans un plat en triangles isoc\`eles est un r\'etract de la famille de tous les triangles de $\G$ le long de $(C_z^r)_{z,r}$.
\end{lemme}

\begin{proof}
Ce r\'esultat est une cons\'equence directe du lemme \ref{reduc}, compte tenu de la d\'efinition de la famille $(C_r^z)_{z,r}$.
\end{proof}

\section{Saturation au sens des plats en triangles isoc\`eles}\label{S-saturation isoc}

Enfin, nous notons que (comme sugg\'er\'e par les lemmes \ref{parallele} (ii) et \ref{reduc}) nous aurions  \'egalement pu ``saturer" l'enveloppe analytique au sens des plats en triangles isoc\`eles. Cette observation permet  de r\'eduire tous les triangles sur des points, m\^eme dans le cas o\`u les points $A,B,C$ sont dans un m\^eme plat en triangles isoc\`eles. Nous d\'ecrivons la r\'eduction correspondante dans le lemme \ref{reduc-sat} ci-dessous. 

Notons \'egalement que le caract\`ere isol\'e  des plats en triangles isoc\`eles n'est pas essentiel pour \'etablir la propri\'et\'e de d\'ecroissance rapide:  m\^eme si les plats en triangles isoc\`eles branchent, on peut quand m\^eme se servir du lemme 15 de \cite{rd}   (plut\^ot que de la proposition 17), pour r\'eduire les produits de convolution \`a des convolutions partielles le long de triangles isoc\`eles, et appliquer alors les techniques de \cite{rd} (adapt\'ees au cas isoc\`ele). 

L'enveloppe analytique satur\'ee est d\'efinie par:

\begin{definition}
Soit $\gamma$ un segment g\'eod\'esique de $X$. On appelle \emph{enveloppe analytique satur\'ee} de $\gamma$, et on note $H'_\sat(\gamma)$, la r\'eunion de l'enveloppe analytique $H(\gamma)$ de $\gamma$ et de tous les plats en triangles isoc\`eles qui ont une intersection non vide avec $H(\gamma)$. 
\end{definition}

Nous tronquons cette enveloppe \'egalement :

\begin{definition}
Soit $\gamma$ un segment g\'eod\'esique de $X$. On appelle \emph{enveloppe analytique satur\'ee r\'eduite} de $\gamma$, et on note $H_\sat'(\gamma)$, l'ensemble des points de l'enveloppe analytique satur\'ee $H_\sat(\gamma)$ de $\gamma$ qui sont \`a distance au plus $2|\gamma|$ du milieu de $\gamma$, o\`u $|\gamma|$ est la longueur CAT(0) de $\gamma$ dans $X$. 
\end{definition}

Ceci fournit une approche alternative pour la propri\'et\'e RD, qui repose sur le crit\`ere de I.\ Chatterji et K.\ Ruane \cite{CR} pour la propri\'et\'e RD, plut\^ot que sur la notion de branchement polynomial pour $G$.  Cette autre approche n'est cependant pas sensiblement diff\'erente \`a celle pr\'esent\'ee  pour  le th\'eor\`eme \ref{wise-branchpol} :  les deux articles \cite{CR} et \cite{rd} pr\'esentent leurs crit\`eres respectifs comme des modifications imm\'ediates de la proposition 2.3  de V.\ Lafforgue \cite{Laf} ; ces deux approches appartiennent  \`a un m\^eme cercle d'id\'ees.

Plus pr\'ecis\'ement, la d\'emonstration du lemme \ref{reduc} s'adapte (imm\'ediatement) au cas satur\'e et implique que:

\begin{lemme}\label{reduc-sat}
Soient $A,B,C$ trois sommets de $X$ et soit $D$ le disque totalement g\'eod\'esique de $X$ de bord le triangle g\'eod\'esique $(ABC)$. Alors l'intersection 
\[
H_\sat'([AB])\cap  H_\sat'([BC])\cap  H_\sat'([AC])
\]
est non vide
\end{lemme}

Ce lemme, combin\'e avec les arguments de la Section \ref{S-deduire RD}, montre alors que la proposition 1.7 de \cite{CR} s'applique, en choisissant pour ensembles $C(x,y)$ l'ensemble des sommets de l'enveloppe analytique satur\'ee r\'eduite du segment  $[x,y]$ de $X$.





\end{document}